\documentclass[aps,epsfig,nofootinbib,article]{revtex4}
\usepackage{mathrsfs}

\usepackage{graphicx}
\usepackage{amsfonts}
\usepackage{epstopdf}
\usepackage{latexsym}
\usepackage{amssymb}
\usepackage{amssymb}
\usepackage{esvect}
\usepackage[center]{subfigure}

\begin{document}
 \newcommand{\bq}{\begin{equation}}
 \newcommand{\eq}{\end{equation}}
 \newcommand{\bqn}{\begin{eqnarray}}
 \newcommand{\eqn}{\end{eqnarray}}
 \newcommand{\nb}{\nonumber}
 \newcommand{\lb}{\label}
\newcommand{\PRL}{Phys. Rev. Lett.}
\newcommand{\PL}{Phys. Lett.}
\newcommand{\PR}{Phys. Rev.}
\newcommand{\CQG}{Class. Quantum Grav.}

\title{A non grid-based interpolation scheme for the eigenvalue problem}

\author{Kai Lin$^{1)}$}\email{lk314159@hotmail.com}
\author{Wei-Liang Qian$^{2,3)}$}\email{wlqian@usp.br}

\affiliation{1) Universidade Federal de Itajub\'a, Instituto de F\'isica e Qu\'imica, Itajub\'a, MG, Brasil}
\affiliation{2) Escola de Engenharia de Lorena, Universidade de S\~ao Paulo, Lorena, SP, Brasil}
\affiliation{3) Faculadade de Engenharia de Guaratinguet\'a, Universidade Estadual Paulista, Guaratinguet\'a, SP, Brasil}

\begin{abstract}
We propose a non grid-based interpolation scheme based on the information from the data collected from the vicinity of the query point.
As a non-grid-based interpolation, the data points can be distributed randomly in a small region, 
and the interpolation is constructed so that it naturally makes use of the information not only on the function value but also on its higher order derivatives. 
The main advantage of the present approach is that the precision of the interpolation can be adjusted in accordance to the quantity of the data, 
in other words, a balance between the precision and the computational cost can be achieved by properly choosing the size of the neighborhood where the data points are collected.
The method is applicable to univariate as well as multivariate functions.
We show that the proposed scheme is efficient and precise.
The present approach is then employed to study the eigenvalue problem.
\\
\\

\end{abstract}

\maketitle
\section{Introduction}
\renewcommand{\theequation}{1.\arabic{equation}} \setcounter{equation}{0}

In physics, eigenvalue problem is associated with many important practical applications.
To name a few, in the classical mechanics, small oscillations of a system about the equilibrium is treated
as superposition of normal modes via the characteristic equation for the frequencies; 
in quantum theory, the atomic spectra can be obtained by solving the time-independent Schrodinger equation;
in general relativity, the stability of a black hole metric is studied by investigating the temporal evolution of small perturbations known as quasinormal modes.
All of the above examples are closely related to the eigenvalue problem.
However, owing to the mathematical difficulties, one usually has to resort to numerical approaches, and therefore, various methods were proposed.

In this work, we propose a non-grid-based interpolation approach.
The interpolation is based on Taylor series, and is constructed by using the information of a set of randomly scattered data points.
We investigate the precision as well as efficiency of the proposed method, by comparing our results to those obtained by standard procedures of {\it MATLAB} or {\it Mathematica}. 
The present interpolation scheme is then employed to discuss the eigenvalue problem.
The paper is organized as follows. 
In section 2, the basic framework of the proposed interpolation method is presented, 
and a few applications are discussed in section 3. 
In section 4, we employ the proposed interpolation scheme to study differential equation and the eigenvalue problem. 
Section 5 is devoted to discussions and conclusions.

\section{A non grid-based interpolation scheme for scattered data points}

\renewcommand{\theequation}{2.\arabic{equation}} \setcounter{equation}{0}

In this section, we propose an interpolation scheme by making use of data points in a small vicinity of the query point. 
Here the data points in question can be scattered, and therefore they do not necessarily sit on the vertices of any regular grid.
The Taylor series \cite{1} for a multivariable function $f(\mathbf{x})\equiv f(x_1,\cdots,x_D)$ about a given point $\mathbf{a}\equiv (a_1,\cdots,a_D)$ read:

 \bqn
 \lb{1}
f(x_1,\cdots,x_D)&=&\sum^\infty_{n_1=0}\cdots\sum^\infty_{n_D=0}\frac{\prod^{D}_{i=1}(x_i-a_i)^{n_i}}{\prod^{D}_{i=1}n_i!}\left(\frac{\partial^{n_1+\cdots+n_D}f}{\partial
x_1^{n_1}\cdots\partial x_D^{n_D}}\right)(a_1,\cdots,a_D)\nb\\
&=&f(a_1,\cdot\cdot\cdot,a_D)+\sum^D_{j=1}\frac{\partial f(a_1,\cdot\cdot\cdot,a_D)}{\partial x_j}(x_j-a_j)\nb\\
&&+\frac{1}{2!}\sum^D_{j=1}\sum^D_{k=1}\frac{\partial^2f(a_1,\cdot\cdot\cdot,a_D)}{\partial x_j\partial x_k}(x_j-a_j)(x_k-a_k)+\cdots
 \eqn
In the special case of a bivariate function, the above series read
 \bqn
 \lb{2}
f(x,y)&=&f(a,b)+(x-a)f_{,x}(a,b)+(y-b)f_{,y}(a,b)\nb\\
&&+\frac{1}{2}(x-a)^2f_{,xx}(a,b)+(x-a)(y-b)f_{,xy}(a,b)+\frac{1}{2}(y-b)^2f_{,yy}(a,b)+\cdots
 \eqn
and for a univariate function, it is simply expressed as
 \bqn
 \lb{3}
f(x)&=&f(a)+(x-a)f'(a)+\frac{1}{2}(x-a)^2f''(a)+\frac{1}{3!}(x-a)^3f'''(a)+\cdots
 \eqn

Since our goal is to carry out an interpolation based on the information of a set of $N$ scattered points distributed in a small neighborhood, let us apply the above Taylor series $N$ times to each one of the data points around query point. 
The results can be written in a form of matrix product as follows
 \bqn
 \lb{4}
F=M D
 \eqn
where $F$ and $D$ are $N\times 1$ column vectors and $M$ is a $N\times N$ matrix. $F$ contains the information from $N$ data points, $D$ contains values of the function $f$ and its derivatives at the query point and the $i$-th line of the matrix $M$ consists of different power functions of the coordinate relative difference between the $i$-th date point and the query point.
For instance, consider a bivariate function $f(x,y)$, and the data points are the function values at coordinates $(x_i,y_i)$ with $i=1,2,\cdots$, $F$, $D$ and $M$ can be written down as follows:

 \bqn
 \lb{5}
F=\left(
    f(x_1,y_1),
    f(x_2,y_2),
    f(x_3,y_3),
    f(x_4,y_4),
    f(x_5,y_5),
    f(x_6,y_6),
    \cdots
\right)^T
 \eqn

 \bqn
 \lb{6}
M= \left(
  \begin{array}{ccccccc}
    1 & x_1-x_0 & y_1-y_0 & \frac{(x_1-x_0)^2}{2} & (x_1-x_0)(y_1-y_0) & \frac{(y_1-y_0)^2}{2} &\cdots \\
    1 & x_2-x_0 & y_2-y_0 & \frac{(x_2-x_0)^2}{2} & (x_2-x_0)(y_2-y_0) & \frac{(y_2-y_0)^2}{2} &\cdots \\
    1 & x_3-x_0 & y_3-y_0 & \frac{(x_3-x_0)^2}{2} & (x_3-x_0)(y_3-y_0) & \frac{(y_3-y_0)^2}{2} &\cdots \\
    1 & x_4-x_0 & y_4-y_0 & \frac{(x_4-x_0)^2}{2} & (x_4-x_0)(y_4-y_0) & \frac{(y_4-y_0)^2}{2} &\cdots \\
    1 & x_5-x_0 & y_5-y_0 & \frac{(x_5-x_0)^2}{2} & (x_5-x_0)(y_5-y_0) & \frac{(y_5-y_0)^2}{2} &\cdots \\
    1 & x_6-x_0 & y_6-y_0 & \frac{(x_6-x_0)^2}{2} & (x_6-x_0)(y_6-y_0) & \frac{(y_6-y_0)^2}{2} &\cdots \\
    \cdots & \cdots & \cdots & \cdots & \cdots & \cdots &\cdots \\
  \end{array}
\right)
 \eqn

 \bqn
 \lb{7}
D= \left(
    f(x_0,y_0),
    f_{,x}(x_0,y_0),
    f_{,y}(x_0,y_0),
    f_{,xx}(x_0,y_0),
    f_{,yx}(x_0,y_0),
    f_{,yy}(x_0,y_0),
    \cdots
\right)^T
 \eqn

In particular, for a univariate function $f(x)$ and $N=4$, one has

 \bqn
 \lb{8}
F=\left(
    f(x_1),
    f(x_2),
    f(x_3),
    f(x_4),
    \cdots
\right)^T
 \eqn

 \bqn
 \lb{9}
M= \left(
  \begin{array}{cccccc}
    1 & x_1-x_0 & \frac{(x_1-x_0)^2}{2} & \frac{(x_1-x_0)^3}{3!} & \frac{(x_1-x_0)^4}{4!} &\cdots \\
    1 & x_2-x_0 & \frac{(x_2-x_0)^2}{2} & \frac{(x_2-x_0)^3}{3!} & \frac{(x_2-x_0)^4}{4!} &\cdots \\
    1 & x_3-x_0 & \frac{(x_3-x_0)^2}{2} & \frac{(x_3-x_0)^3}{3!} & \frac{(x_3-x_0)^4}{4!} &\cdots \\
    1 & x_4-x_0 & \frac{(x_4-x_0)^2}{2} & \frac{(x_4-x_0)^3}{3!} & \frac{(x_4-x_0)^4}{4!} &\cdots \\
    \cdots & \cdots & \cdots & \cdots & \cdots &\cdots \\
  \end{array}
\right)
 \eqn

 \bqn
 \lb{10}
D= \left(
    f(x_0),
    f'(x_0),
    f''(x_0),
    f'''(x_0),
    \cdots
\right)^T
 \eqn

This implies that the column vector $D$ can be expressed in terms of $F$ and $M$ once $M$ has a nonzero determinant, $D=A^{-1}F$. In particular, if only a few particular elements of $D$ is required, one may use Cramer's rule to evaluate them: $D_i=\det(M_i)/\det(M)$ where $M_i$ is the matrix formed by replacing the $i$-th column of $M$ by the column vector $A$. For instance, in the above example, $f''(x_0)=\det(M_3)/\det(M)$. It is noted that if the (partial) derivatives instead of the function value are known, Eq.(\ref{4}) and the subsequent procedure can be carried out in a very similar fashion.

\section{Precision and efficiency of the method}
\renewcommand{\theequation}{3.\arabic{equation}} \setcounter{equation}{0}

\begin{figure*}
\includegraphics[width=8cm]{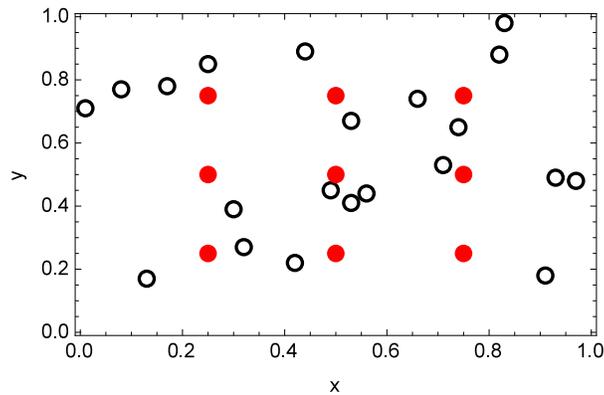}
\caption{The black empty circles are the 21 random points in Table \ref{TableI},
and the filled red circles correspond the points in Table \ref{TableII}.} \label{fig3}
\end{figure*}

Because the information on all the data points are expressed in terms of Taylor series by Eq.(\ref{4}), it is expected that the present interpolation scheme inherits the efficiency of any method which is also based on Taylor series.
In what follows, we study the precision of the method by applying it to few examples.
First, let us consider a bivariate function $f(x,y)=\sin(xy^2)$. 
We randomly sample 21 points in the region $0\le x\le1$ and $0\le y\le1$, and use the function value at those points to calculate all first order partial derivatives of $f(x,y)$ using the proposed interpolation scheme.
The results are shown in Table \ref{TableI} below:
\begin{table}[ht]
\caption{\label{TableI} A comparison of the calculated $\frac{\partial f}{\partial x}$ and $\frac{\partial f}{\partial y}$ with those from the analytic expression for the 21 random points.}
\begin{tabular}{cccccccc}
         \hline
$(x,y)$ & $f=\sin(xy^2)$ & $\bar{f}_x$ &$\tilde{f}_x$&$\Delta f_x$ &
$\bar{f}_y$& $\tilde{f}_y$&$\Delta f_y$
        \\
        \hline
 \{0.44,0.89\} & 0.341511 & 0.744477 & 0.744112 & -0.000490826 & 0.736112
   & 0.718355 & -0.0241236 \\
 \{0.83,0.98\} & 0.715355 & 0.67109 & 0.712135 & 0.0611608 & 1.13674 &
   1.13167 & -0.00446141 \\
 \{0.3,0.39\} & 0.0456142 & 0.151942 & 0.153233 & 0.00849649 & 0.233756 &
   0.232836 & -0.00393799 \\
 \{0.53,0.67\} & 0.235679 & 0.436255 & 0.436543 & 0.000659489 & 0.690194
   & 0.692824 & 0.00381002 \\
 \{0.01,0.71\} & 0.00504098 & 0.504094 & 0.501928 & -0.00429524 &
   0.0141998 & 0.0204939 & 0.443253 \\
 \{0.32,0.27\} & 0.0233259 & 0.0728802 & 0.0727901 & -0.00123588 &
   0.172753 & 0.172403 & -0.00202474 \\
 \{0.97,0.48\} & 0.221632 & 0.22467 & 0.229766 & 0.0226816 & 0.908041 &
   0.928144 & 0.0221382 \\
 \{0.42,0.22\} & 0.0203266 & 0.04839 & 0.0473943 & -0.020577 & 0.184762 &
   0.180365 & -0.0237969 \\
 \{0.08,0.77\} & 0.0474142 & 0.592233 & 0.590927 & -0.00220542 & 0.123061
   & 0.122278 & -0.00636394 \\
 \{0.53,0.41\} & 0.0889752 & 0.167433 & 0.166948 & -0.00289553 & 0.432876
   & 0.433456 & 0.00133978 \\
 \{0.91,0.18\} & 0.0294797 & 0.0323859 & 0.0179674 & -0.445211 & 0.327458
   & 0.277667 & -0.152052 \\
 \{0.17,0.78\} & 0.103244 & 0.605149 & 0.606371 & 0.00201958 & 0.263783 &
   0.264328 & 0.00206698 \\
 \{0.25,0.85\} & 0.179644 & 0.710746 & 0.715137 & 0.00617731 & 0.418086 &
   0.409043 & -0.0216305 \\
 \{0.82,0.88\} & 0.593184 & 0.623444 & 0.637844 & 0.0230979 & 1.16187 &
   1.16065 & -0.00105318 \\
 \{0.13,0.17\} & 0.00375699 & 0.0288998 & 0.0430017 & 0.487959 &
   0.0441997 & 0.0190376 & -0.569282 \\
 \{0.66,0.74\} & 0.353599 & 0.512223 & 0.510209 & -0.003932 & 0.913696 &
   0.914023 & 0.000357864 \\
 \{0.74,0.65\} & 0.307581 & 0.402018 & 0.403238 & 0.00303393 & 0.915364 &
   0.91452 & -0.000922349 \\
 \{0.56,0.44\} & 0.108204 & 0.192463 & 0.192158 & -0.00158502 & 0.489907
   & 0.490137 & 0.000469295 \\
 \{0.71,0.53\} & 0.198119 & 0.275332 & 0.275418 & 0.000312201 & 0.737682
   & 0.738331 & 0.000879517 \\
 \{0.49,0.45\} & 0.0990623 & 0.201504 & 0.201672 & 0.000835784 & 0.438831
   & 0.43798 & -0.00193901 \\
 \{0.93,0.49\} & 0.221442 & 0.234139 & 0.23768 & 0.0151232 & 0.888773 &
   0.903173 & 0.0162025 \\
        \hline
\end{tabular}
\end{table}
where $\tilde{f}_x$ and $\tilde{f}_y$ are the calculated values of $\frac{\partial f(x,y)}{\partial x}$ and $\frac{\partial f(x,y)}{\partial y}$ by using the proposed interpolation scheme,
$\bar{f}_x$ and $\bar{f}_y$ are those obtained by the analytic expressions, 
and $\Delta{f}_x=\frac{\tilde{f}_x-\bar{f}_x}{\bar{f}_x}$ $\Delta{f}_y=\frac{\tilde{f}_y-\bar{f}_y}{\bar{f}_y}$ are the relative difference between the exact and the numerical values.
Since we have only made use of 21 points, one can only expand the function up to its fourth order partial derivatives. 
However, we find that the resultant precision is reasonable, provided that these few sampled points were obtained randomly.

In the second example, we apply the method to further evaluate the function and its partial derivatives at any given point. 
We made use of the same information of the above 21 point, and evaluate the function value as well as the derivative $\frac{\partial^2 f}{\partial x\partial y}$ at the rectangular grid points shown in Fig.\ref{fig3},
the results are presented in Table \ref{TableII}.
\begin{table}[ht]
\caption{\label{TableII} A comparison of the calculated $f$ and $\frac{\partial^2 f}{\partial x\partial y}$ with those from the analytic expression at 9 points.}
\begin{tabular}{ccccccc}
         \hline
$(x,y)$ &$\bar{f}$ &$\tilde{f}$&$\Delta f$ & $\bar{f}_{xy}$&
$\tilde{f}_{xy}$&$\Delta f_{xy}$
        \\
        \hline
\{0.25,0.25\} & 0.0156244 & 0.0155806 & -0.00279899 & 0.499817 &
0.517806
   & 0.035991 \\
 \{0.25,0.5\} & 0.0624593 & 0.062147 & -0.00500034 & 0.994144 & 0.982991 &
   -0.0112189 \\
 \{0.25,0.75\} & 0.140162 & 0.140196 & 0.00024389 & 1.45563 & 1.47183 &
   0.0111322 \\
 \{0.5,0.25\} & 0.0312449 & 0.0311037 & -0.00451952 & 0.499268 & 0.483332
   & -0.0319184 \\
 \{0.5,0.5\} & 0.124675 & 0.124621 & -0.000432479 & 0.976613 & 0.986277 &
   0.00989526 \\
 \{0.5,0.75\} & 0.277557 & 0.277861 & 0.00109612 & 1.32397 & 1.30442 &
   -0.0147678 \\
 \{0.75,0.25\} & 0.0468578 & 0.0463861 & -0.0100679 & 0.498353 & 0.461077
   & -0.0747966 \\
 \{0.75,0.5\} & 0.186403 & 0.186352 & -0.000275179 & 0.947523 & 0.977499 &
   0.0316367 \\
 \{0.75,0.75\} & 0.409472 & 0.409368 & -0.000253452 & 1.10937 & 1.10065 &
   -0.0078521 \\
        \hline
\end{tabular}
\end{table}
where $\tilde{f}$ and $\tilde{f}_{xy}$ represent the calculated values of the function and $\frac{\partial^2f}{\partial x\partial y}$ by using the proposed interpolation scheme, 
$\bar{f}$ and $\bar{f}_{xy}$ are those obtained by the analytic expressions,
and $\Delta f$ $\Delta f_{xy}$ are the relative difference between the exact and the numerical values.

In the following example, we employ the interpolation method to find the extrema of an unknown function. 
This can be achieved by solving the linear equation where the interpolated first order derivative equals zero.
Let us consider the Jacobi elliptic function, $\text{sn}(x|\frac{1}{3})$.
Our goal is to find the extreme between $(1,7)$ by using the function values at 7 points distributed evenly between the above interval.
The function possesses two extrema in the interval, a maximum at $5.20175$ and a minimum at $1.73392$.
Comparing to the results obtained by evoking the ``Polynomial" command of {\it Mathematica}, $5.1633$ and $1.9322$, the relative differences are $-0.00739242$和$0.114378$.
Our interpolation scheme gives $5.15693$ and $1.66738$, the relative differences are $-0.00900098$ and $-0.0383723$, which are reasonably good.

As another example, we show that the proposed method also can be carried out when the values of the derivative of an unknown function are given at some points.
Let us assume that the function $f(x)=(x-1)\sin(x)\sin(x^{-2})$ is unknown, but its values at $x=1,2,3,4,5,6,7,8,9,10$ are given. 
In addition, it is provided that the extreme are at the following three points $x_0=1.77251,4.55625,7.74447$.
Our goal is to evaluate the values of the function at the above three points.
In this case, we are provided by some additional information on the derivative of the function.
As discussed by the end of section I, we can make use of these information by adding three more elements, namely, three zeros correspond to the first order derivative at those extrema, to the column vector $F$ in Eq.(\ref{4}).
The right hand side of Eq.(\ref{4}) can be subsequently written down by expanding the first order derivative in terms of higher order derivatives.
The function values according to the analytic expression are $f(x_0)=0.23685,-0.169158,0.111772$. 
In comparison to those obtained by the ``InterpolatingPolynomial" command of {\it Mathematica}, $f(x_0)=0.235812,-0.169107,0.111709$, the relative differences are $-0.00438185,-0.000301373,-0.000571772$
The calculated function values by the proposed method are $f(x_0)=0.236839,-0.169167,0.111774$, the relative differences are $-0.000045039,-0.0000482994,-0.0000180876$.
Again, the results turn out to be of the same precision of those obtained by {\it Mathematica}.

\section{Applications in differential equation and eigenvalue problem}
\renewcommand{\theequation}{4.\arabic{equation}} \setcounter{equation}{0}

In numerical computation, many methods have been proposed to solve differential equation.
To this day, to further improve the precision and the efficiency in solving differential equation still remains as a challenging as well as compelling topic to the community.
Differential quadrature method (DQM) \cite{3} was proposed in the 70s by Bellman and Casti, it is used as one of the widely recognized methods for numerical differentiation.
The method proposed here is able to achieve the same precision and efficiency of the former, as shown below.

To apply the proposed method, first let us revisit the Taylor series of $f_i(x_1,\cdots,x_N)$ around $(a_1,\cdots,a_N)$ Eq.(\ref{1}):

 \bqn
 \lb{11}
f_i(x_1,\cdots,x_D)&=&f(a_1,\cdots,a_D)+\sum^D_{j=1}\frac{\partial f(a_1,\cdots,a_D)}{\partial x_j}(x_j-a_j)\nb\\
&&+\frac{1}{2!}\sum^D_{j=1}\sum^D_{k=1}\frac{\partial^2f(a_1,\cdots,a_D)}{\partial x_j\partial x_k}(x_j-a_j)(x_k-a_k)+\cdots \nb
 \eqn
 
For a univariate function, one has \cite{4}：

 \bqn
 \lb{12}
\left(
\begin{array}{cccccc}
    f(x_1)\\
    f(x_2)\\
    f(x_3)\\
    f(x_4)\\
    \cdots\\
    \end{array}
\right)=
\left(
\begin{array}{cccccc}
    f(x_0)\\
    f(x_0)\\
    f(x_0)\\
    f(x_0)\\
    \cdots\\
    \end{array}
\right)+ 
\left(
  \begin{array}{cccccc}
     x_1-x_0 & \frac{(x_1-x_0)^2}{2} & \frac{(x_1-x_0)^3}{3!}& \frac{(x_1-x_0)^4}{4!} & \cdots \\
     x_2-x_0 & \frac{(x_2-x_0)^2}{2} & \frac{(x_2-x_0)^3}{3!}& \frac{(x_2-x_0)^4}{4!} & \cdots \\
     x_3-x_0 & \frac{(x_3-x_0)^2}{2} & \frac{(x_3-x_0)^3}{3!}& \frac{(x_3-x_0)^4}{4!} & \cdots \\
     x_4-x_0 & \frac{(x_4-x_0)^2}{2} & \frac{(x_4-x_0)^3}{3!}& \frac{(x_4-x_0)^4}{4!} & \cdots \\
     \cdots & \cdots & \cdots & \cdots & \cdots \\
  \end{array}
\right) \left(
\begin{array}{cccccc}
    f'(x_0)\\
    f''(x_0)\\
    f'''(x_0)\\
    f''''(x_0)\\
    \cdots \\
    \end{array}
\right)
 \eqn

Now, let us interpret the above expression from its right hand side to the left hand side, in other words, to view the function $f(x)$ as well as its derivatives at a given point $x_0$ as linear combination of function values at $(x_1,\cdots,x_N)$.
In order to solve the differential equation, we adopt the assumption that the Taylor series expansion is valid in the domain of interest\footnote{If not, coordinate transformation usually can be carried out to achieve this requirement.}. 
Then, one divides the domain of the function by using regular grids, or when it is necessary, scattered points. 
These grid points as well as the boundary are associated with the above-mentioned points $(x_1,\cdots,x_N)$ on the right hand side of Eq.(\ref{12}). 
By substituting the left hand side of Eq.(\ref{12}) into the differential equation, one finally obtains a linear equation in terms of the function values at $(x_1,\cdots,x_N)$.
Since the resulting equation is valid for any point in the domain, we can choose $N$ different points for $x_0$ at will and solve the problem by methods of linear algebra.

As an example, let us solve the differential equation $f''(x)+\sin^3(x)=0$ with the boundary conditions $f(0)=3$ and $f(1)=2$. In this example, we discretize the interval $[0,1]$ by using a total of 14 uniformly distributed points, and solve the corresponding problem in linear algebra involving 14 linear equations. The results are shown in Table \ref{TableIII} where exact solution $\bar{f}$ is compared to the numerical solution $\tilde{f}$, the proposed method has achieved a very high precision.
\begin{table}[ht]
\caption{\label{TableIII} A comparison between the exact solution $\bar{f}$ and numerical solution $\tilde{f}$ of the differential equation $f''(x)+\sin^3(x)=0$ with $f(0)=3$ and $f(1)=2$, $\Delta f$ is the relative differences.}
\begin{tabular}{cccccccc}
         \hline
$x$ & $\bar{f}$ &$\tilde{f}$&$\Delta f$
        \\
        \hline
 0 & 3. & 3. & 0. \\
 $\frac{1}{13}$ & 2.92611 & 2.92611 & 6.213364347442804$\times10^{-13}$ \\
 $\frac{2}{13}$ & 2.85222 & 2.85222 & 6.21084423474243$\times10^{-13}$ \\
 $\frac{3}{13}$ & 2.77831 & 2.77831 & 6.404848458855489$\times10^{-13}$ \\
 $\frac{4}{13}$ & 2.70432 & 2.70432 & 6.583365058193686$\times10^{-13}$ \\
 $\frac{5}{13}$ & 2.63016 & 2.63016 & 6.777424923871763$\times10^{-13}$ \\
 $\frac{6}{13}$ & 2.55569 & 2.55569 & 6.984482190814282$\times10^{-13}$ \\
 $\frac{7}{13}$ & 2.48069 & 2.48069 & 7.205502781153421$\times10^{-13}$ \\
 $\frac{8}{13}$ & 2.40488 & 2.40488 & 7.440020803985056$\times10^{-13}$ \\
 $\frac{9}{13}$ & 2.32793 & 2.32793 & 7.691674102109481$\times10^{-13}$ \\
 $\frac{10}{13}$ & 2.24944 & 2.24944 & 7.976359101089474$\times10^{-13}$ \\
 $\frac{11}{13}$ & 2.16895 & 2.16895 & 8.25879045958274$\times10^{-13}$ \\
 $\frac{12}{13}$ & 2.08598 & 2.08598 & 8.827602433319681$\times10^{-13}$ \\
 1 & 2. & 2. & 0. \\
        \hline
\end{tabular}
\end{table}

Now we are in a position to apply the proposed method is to the eigenvalue problem. 
Without loss of generality, we present the procedure by solving the stationary solution of one-dimensional Schrodinger equation.
As in the above scheme for solving differential equation, by properly discretizing the domain of the unknown function, one obtains an array of linear equations in terms of the function values at $N$ grid points.
The key characteristic of the eigenvalue problem is that resultant linear equations are homogeneous, so that it can be written into the form
\bqn
M(\omega) \phi = 0
\eqn
where $\phi$ is a $N \times 1$ column vector consists of the function values at $N$ grid points, and $M$ is a $N \times N$ matrix and $\omega$ is the eigenvalue. 
This implies that $\det(M(\omega))=0$, where $\det(M(\omega))$ is a polynomial of $\omega$ with known coefficients.

As a first example let us consider the one dimensional infinity potential well. The corresponding eigenvalue equation reads
 \bqn
 \lb{13}
&&\phi''(x)+\omega^2\phi(x)=0,~~~0\le x\le 1\nb\\
&&\phi(1)=\phi(0)=0
 \eqn
where the eigenvalues are known to be $\omega_n=n\pi$. Let us again consider a total of 14 uniformly distributed grid points in $[0,1]$ and then solve the polynomial equation for $\omega$ numerically.
We obtain the following real eigenvalues $3.14159, 6.28319, 9.42387, 12.5465$. We conclude that the resulting precision is reasonably good and it is higher with smaller $n$.

In the second example, we study the Schrodinger equation with the Poschl-Teller potential defined as follows:
 \bqn
 \lb{14}
U(x)&=&\frac{V_0}{2}\left[\frac{k(k-1)}{\sin^2(\alpha
x)}+\frac{\lambda(\lambda-1)}{\cos^2(\alpha x)}\right],\nb\\
V_0&=&\frac{\hbar^2\alpha^2}{m_0}
 \eqn
where the eigenvalues are known to be $E_n=\frac{V_0}{2}\left(k+\lambda+2n\right)$.
Here we assume $\alpha=\frac{\pi}{2}$, and solve the eigenvalue problem with the proposed method, and compare the exact values $\bar{E}_n$ with the numerical results $\tilde{E}_n$ in Table \ref{TableIV}.
\begin{table}[ht]
\caption{\label{TableIV} A comparison of the eigenvalues for the Poschl-Teller potential between the exact values $\bar{E}_n$ and the numerical results $\tilde{E}_n$.}
\begin{tabular}{cccccccc}
         \hline
$(k,\lambda)$ &
$\left\{\frac{\bar{E}_1}{V_0},\frac{\tilde{E}_1}{V_0}\right\}$ &
$\left\{\frac{\bar{E}_2}{V_0},\frac{\tilde{E}_2}{V_0}\right\}$&$\left\{\frac{\bar{E}_3}{V_0},\frac{\tilde{E}_3}{V_0}\right\}$
        \\
        \hline
$(2,2)$&$\left\{18,~~18.0492\right\}$&$\left\{32,~~31.8772\right\}$&$\left\{50,~~49.4934\right\}$\\
$(3,2)$&$\left\{24.5,~~24.4631\right\}$&$\left\{40.5,~~40.546\right\}$&$\left\{60.5,~~60.2104\right\}$\\
$(2,3)$&$\left\{24.5,~~24.4974\right\}$&$\left\{40.5,~~40.4348\right\}$&$\left\{60.5,~~60.45\right\}$\\
$(3,3)$&$\left\{32,~~31.9963\right\}$&$\left\{50,~~49.9572\right\}$&$\left\{72,~~71.2222\right\}$\\
$(10,10)$&$\left\{242,~~241.792\right\}$&$\left\{288,~~289.955\right\}$&$\left\{338,~~330.614\right\}$\\
        \hline
\end{tabular}
\end{table}
We see that reasonable precision is obtained.

Let us consider another example of the following trigonometric potential
 \bqn
 \lb{15}
U(x)=V_a\cot\left(\frac{\pi}{\alpha}x\right)
 \eqn
with the following known analytic eigenvalues
 \bqn
 \lb{15}
E_n&=&\frac{\pi^2}{A\alpha^2}(n^2+4n\lambda_a-2\lambda_a),\nb\\
\lambda_a&=&\frac{1}{4}\left(\sqrt{\frac{4AV_a\alpha^2}{\pi^2}+1}-1\right)
 \eqn
where
$A=\frac{2m_0}{\hbar^2}$. 
For simplicity, we choose $\alpha=A=1$, and compare the exact values $\bar{E}_n$ with the numerical results $\tilde{E}_n$ in Table \ref{TableV}.
\begin{table}[ht]
\caption{\label{TableV} A comparison of the eigenvalues for the trigonometric potential between the exact values $\bar{E}_n$ and the numerical results $\tilde{E}_n$.}
\begin{tabular}{cccccccc}
         \hline
$V_a$ & $\left\{\bar{E}_1,\tilde{E}_1\right\}$ &
$\left\{\bar{E}_2,\tilde{E}_2\right\}$&
$\left\{\bar{E}_3,\tilde{E}_3\right\}$&
$\left\{\bar{E}_4,\tilde{E}_4\right\}$
        \\
        \hline
$1$&$\left\{10.7847,~~11.3125\right\}$&$\left\{42.2239,~~40.3267\right\}$&$\left\{93.4022,~~95.4006\right\}$&$\left\{164.32,~~160.329\right\}$\\
$10$&$\left\{16.0275,~~15.9083\right\}$&$\left\{57.9522,~~58.0788\right\}$&$\left\{119.616,~~119.774\right\}$&$\left\{201.019,~~196.455\right\}$\\
$50$&$\left\{27.6907,~~27.6959\right\}$&$\left\{92.9418,~~92.8824\right\}$&$\left\{177.932,~~177.894\right\}$&$\left\{282.662,~~278.542\right\}$\\
$100$&$\left\{36.7359,~~36.7354\right\}$&$\left\{120.077,~~120.087\right\}$&$\left\{223.158,~~222.993\right\}$&$\left\{345.978,~~343.792\right\}$\\
$1000$&$\left\{104.403,~~104.402\right\}$&$\left\{323.079,~~323.075\right\}$&$\left\{561.494,~~561.621\right\}$&$\left\{819.649,~~820.084\right\}$\\
        \hline
\end{tabular}
\end{table}
We again arrive the conclusion that the proposed method yields satisfactory results.

\section{Discussions and conclusions}
\renewcommand{\theequation}{5.\arabic{equation}} \setcounter{equation}{0}

In this work, we propose a non-grid-based interpolation scheme based on the information on the data in the vicinity of the query point.
As a non-grid-based interpolation, the data points do not have to be distributed uniformly in the sampling area.
Therefore, a very important feature of the proposed method is that the interpolation is not restricted to the data sampling process and 
as a matter of fact, the precision of the interpolation can be adjusted in accordance to the quantity of the data.
We applied the method to differential equation as well as to the eigenvalue problem.

There are many other applications of the eigenvalue problem. 
As discussed above, the study of small perturbations of black hole, known as quasinormal mode, is one of such problems. 
And due to its implication in general relativity and particle physics, it has aroused increasing interest of physicists in recent years.
In this case, the eigenvalue is a complex number therefore the above method cannot be applied straightforwardly.
Other applications include geometric measurements. 
In this context, due to practical restrictions, the proposed method may become quite advantageous since no rectangular grids are necessary.
As a result, measurements can be carried out with much convenience as well as better precision.
Such implementation is desirable and will be carried out in the near future.

\section*{\bf Acknowledgements}

This work is supported in part by Brazilian funding agencies FAPESP, FAPEMIG, CNPq, CAPES, and by Chinese funding agencies NNSFC.


\end{document}